\documentclass[reqno]{amsart}
\usepackage{amsmath,amsfonts,epsfig,amssymb,graphics}

{\theoremstyle{plain}
\newtheorem{thm}{Theorem}[section]
\newtheorem{cor}[thm]{Corollary}
\newtheorem{prop}[thm]{Proposition}
}
{\theoremstyle{definition}
\newtheorem{defn}[thm]{Definition}
\newtheorem{rem}[thm]{Remark}
}

\newcommand{\CC}{{\mathcal C}}
\newcommand{\Sp}{{\mathbb S}}
\newcommand{\G}{{\mathcal G}}
\newcommand{\PP}{{\mathcal P}}
\newcommand{\arcsinh}{\mathrm{arcsinh}}
\newcommand{\arccosh}{\mathrm{arccosh}}

\begin{document}

\title{Collars and partitions of hyperbolic cone-surfaces \ }

\author{Emily B. Dryden}
\address{Department of Mathematics, Bucknell University, Lewisburg, PA 17837} 
\email{ed012@bucknell.edu}
\author{Hugo Parlier}
\address{Section de Math\'ematiques, Universit\'e de Gen\`eve, Switzerland}
\email{hugo.parlier@math.unige.ch}
\thanks{The first author was partially supported by the U.S. National
  Science Foundation grant DMS-0306752.  The second author was
  supported by the Swiss National Science Foundation grants 21 -
  57251.99 and 20 - 68181.02}


\begin{abstract}
For compact Riemann surfaces, the collar theorem and Bers'
partition theorem are major tools for working with simple closed
geodesics.  The main goal of this paper is to prove similar
theorems for hyperbolic cone-surfaces.  Hyperbolic two-dimensional orbifolds
are a particular case of such surfaces.  We
consider all cone angles to be strictly less than $\pi$ to be
able to consider partitions. 
\end{abstract}


\maketitle

\section{Introduction} \label{Sect:S1}
Compact hyperbolic cone-surfaces are a natural
generalization of both compact hyperbolic Riemann surfaces and
compact hyperbolic two-dimensional orbifolds. Both the collar
theorem and Bers' theorem for partitions are very important tools
for Riemann surfaces. It is thus natural to try to find equivalent
theorems for two-dimensional orbifolds, and more generally for cone-surfaces. 

The collar theorem for compact hyperbolic Riemann surfaces (e.g.
\cite{kee74}, \cite{bu78}) states that surrounding a simple closed
geodesic there is a tubular neighborhood, called a collar, which
is a topological cylinder. This neighborhood is of a certain
width which depends uniquely on the length of the geodesic.
Furthermore, if two simple closed geodesics do not
intersect then their collars are disjoint. Finally, the values
given for the widths of the collars are optimal (\cite{pa041}, \cite{sesobook}). 
There has been interest in
proving a similar theorem for orbifolds (e.g. \cite{diphd},
\cite{gema97} and \cite{mat76}), where the object was often to
estimate minimal distance between singular points based on the
order of the points. The collar theorem for hyperbolic
cone-surfaces has the same properties as
the original collar theorem, and is thus a natural generalization.  
Specifically, given a collection of disjoint simple closed geodesics
on a hyperbolic cone-surface $M$, we can complete it to a
partition of $M$, or a set of pairwise disjoint simple closed
geodesics which divide the surface into pairs of pants.
We find collars surrounding these simple closed geodesics
and the cone points in $M$; the collars are disjoint.  
 The width of the collar about a simple closed
geodesic depends on the length of the
geodesic and the size of the largest cone angle in the surface, while
the width of the collar about a cone point depends on the cone angle
at that point.  The collars about geodesics are topological cylinders,
and those about cone points are topological cones.  The values given
for the widths of these collars are optimal.  
The techniques used in the proof are based on hyperbolic
trigonometry and both geometric and topological
properties of simple closed geodesics.

Bers' partition theorem states that, on a compact hyperbolic
Riemann surface, a partition of the surface into pairs of pants
can be chosen such that the lengths of the partitioning geodesics
are bounded by a constant depending solely on the genus. The
optimal constant is not known, although an upper and lower bound
are known (e.g. \cite{bubook}) and certain properties of surfaces
that would realize Bers' constant are also known (e.g.
\cite{pa044}). We find a constant analogous to Bers'
partition constant that depends on the genus and the number of
singular points; any such bound will necessarily depend upon the
number of singular points. The techniques used to prove the
theorem include area arguments using polar coordinates.  An explicit
bound is found in the proof.

In the setting of compact hyperbolic Riemann surfaces, Bers'
theorem has proved to be a useful tool in studying spectral
questions. In particular, it has been used to find a rough
fundamental domain for the action of the Teichm\"{u}ller modular
group, to find an explicit bound on the size of isospectral
families, and in estimates involving Fenchel-Nielsen parameters
(see \cite{bubook}). Its utility stems from the fact that it
allows one to significantly restrict the allowed lengths of
partitioning geodesics. 

The paper is organized as follows.  We begin by giving the
necessary background on hyperbolic cone-surfaces.
In section 3, we prove the collar theorem (Theorem
\ref{thm:collartheorem}) and discuss its consequences.  Bers' theorem (Theorem
\ref{thm:bersthm}) is proved in section 4. \\

The authors are grateful to  
Peter Buser and Carolyn Gordon for reading and providing feedback on  
early drafts of this work.  Much of this research was done while the first  
named author was visiting EPFL, and she thanks the geometry group  
there for their hospitality.

\section{Preliminaries} \label{Sect:S2}

Although most results, definitions and properties in this section
are readily accepted and known, we have not found all of them in
the literature. For this reason, we have tried to be rather
complete.

\begin{defn}
A \emph{hyperbolic cone-surface} is a two-dimensional manifold, $M$,
which can be triangulated by hyperbolic geodesic triangles.
\end{defn}

We are interested in the compact and orientable case and
henceforth will always assume our manifolds to be compact and
orientable. Such manifolds can have a certain number of singular
points and the set of all such points will be denoted by
$\Sigma$. If $\Sigma=\emptyset$ then $M$ is a compact hyperbolic
Riemann surface. If $\Sigma\neq\emptyset$ then $M\setminus
\Sigma$ has a smooth (although not complete) Riemannian metric
with constant curvature $-1$. At each singular point $p$, there is a
collection of triangles having $p$ as a vertex and with given
interior angles at $p$.  
The sum of these given angles is the cone angle associated to $p$.  
In what follows, these angles will be
considered strictly inferior to $\pi$ for two reasons. First, two
singular points with cone angle $\pi$ can be arbitrarily close to
each other or to a simple closed geodesic. Thus there can be no
natural generalization of the collar theorem from compact
hyperbolic Riemann surfaces to hyperbolic cone-surfaces with such singular
points. Furthermore, the notion of partition changes
significantly when we allow cone angle $\pi$, making it difficult
to find a natural generalization of Bers' theorem. If cone angles
are all of the form $2\pi/k$ where $k$ is an integer, then $M$ is
an orbifold. Finally, if one lets the cone angles tend towards
$0$, cusps appear. The case of cusps is well known (\cite
{bubook}) and thus is not mentioned in the theorem. However, our
techniques apply to this case as well.

Our objects of study will be called {\it admissible
cone-surfaces}. These are compact orientable hyperbolic
cone-surfaces which satisfy certain conditions.
First, all cone angles are less than $\pi$. We also require that
$(g,n) \geq (0,4)$, where $g$ is the genus of the surface and $n$
is the number of cone points. Here the ordering is lexicographic
ordering; that is, $(g,n) \geq (0,4)$ means that $ g \geq 0$, and
if $g=0$, then $n \geq 4$. Note that $(g,n)$ cannot equal
$(1,0)$, as we require our surfaces to be endowed with a
hyperbolic metric.

Two curves are considered freely homotopic if they are freely
homotopic on $M \setminus \Sigma$. A closed curve and a cone point $p$
are freely homotopic if they are freely homotopic on $( M \setminus
\Sigma ) \cup \{ p\}$.
Furthermore, by curve or geodesic we mean the set of all points
lying on the curve. The following properties for simple closed
curves are a result of the constant negative curvature, and can
be proved in an identical fashion to their equivalent for Riemann
surfaces.

\begin{prop}
Let $M$ be an admissible cone-surface.
\begin{enumerate}
\item
Every non-trivial simple closed curve on $M\setminus \Sigma$ is freely
homotopic to either a unique simple closed geodesic or a unique cone
point. (For a closed curve $\delta$, the associated closed geodesic will be
denoted $\G(\delta)$.)
\item
If two non-trivial closed curves $\alpha$ and $\beta$ intersect $n$ times,
then $\G(\alpha)$ and $\G(\beta)$ intersect at most $n$ times, or
$\G(\alpha)=\G(\beta)$.
\item
Given two non-intersecting smooth simple curves $\alpha$ and
$\beta$ on $M$ there is at least one geodesic path $c$ between
them such that $d_M(\alpha,\beta)$ is realized by $c$. Such a path
$c$ is perpendicular to $\alpha$ and $\beta$. If $\alpha$ and
$\beta$ are geodesic, in a free homotopy class of paths with end
points moving on $\alpha$ and $\beta$, such a path $c$ is unique.
This property remains true for singular points in place of one or
both geodesics. \end{enumerate}
\end{prop}

\begin{figure}
\center
\epsfxsize=60mm 
\epsfbox{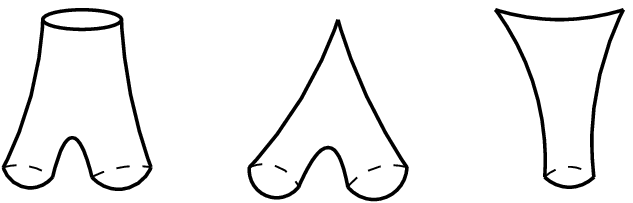}
\caption{}
\label{fig:pairsofpants}
\end{figure}

\begin{figure}
\center
\epsfxsize=40mm 
\epsfbox{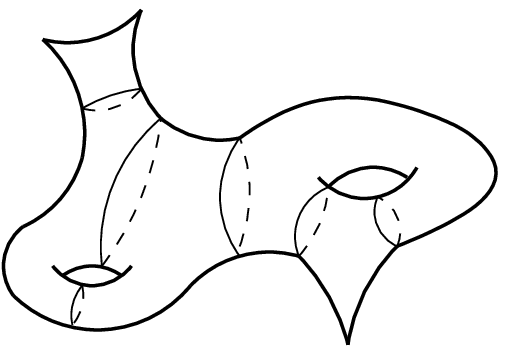}
\caption{}
\label{fig:decompexample}
\end{figure}

Based on properties in the previous proposition, one can define a
partition to be a set of simple closed geodesics
$\PP=\{\gamma_1,\ldots,\gamma_m\}$ such that $M\setminus\PP$ is a
set of pairs of pants. A pair of pants is defined as for smooth
surfaces, except a boundary geodesic can be replaced by a cone
point. We have $m=3g-3+n$ where $g$ is the genus and $n$ the
number of singular points. If a pair of pants has three geodesics
as boundary, then it will be called a $Y$-piece. 
If it has one cone point and two geodesics as boundary then it will be called a
$V$-piece. 
Finally, if the pair of pants has two cone points and one
geodesic as boundary then it will be called a {\it joker's hat}. 
See Figure \ref{fig:pairsofpants} for examples of the various types of pairs of
pants, and Figure \ref{fig:decompexample} for an example of a partition.

The only situation in which a pair of pants has three cone points
as boundary is if $M$ is a sphere with three cone points. For
other basic properties of hyperbolic cone-surfaces see
\cite{cohokebook}.

It is well known that a $Y$-piece can be decomposed into two
isometric right-angled hexagons in the hyperbolic plane. In the
same fashion, a $V$-piece can be decomposed into two isometric
pentagons with four right angles (\cite{diphd}), and a joker's hat
into two isometric quadrilaterals with two right angles. Both such
hyperbolic pentagons and quadrilaterals can in turn be broken
down into two (not necessarily isometric) trirectangles, which are
quadrilaterals with three right angles. For hyperbolic formulas
for these polygons, see \cite{bubook} or \cite{fenibook}.

\section{The collar theorem} \label{Sect:S3}

\begin{figure}
\center
\epsfxsize=30mm 
\epsfbox{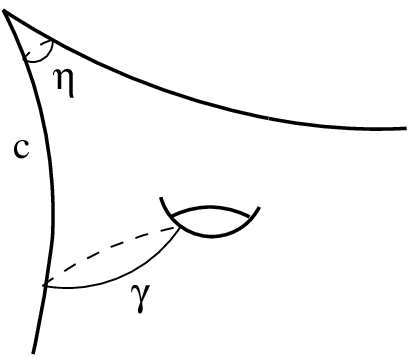}
\caption{}
\label{fig:uniquepants}
\end{figure}

Let $p$ be a cone point and $\gamma$ a simple closed geodesic on
$M$. Let $c$ be a simple geodesic path from $p$ to $\gamma$,
perpendicular to $\gamma$. These elements describe a unique pair
of pants in the following manner. Let $\delta$ be the closed
curve obtained by taking $\gamma\circ c \circ \eta \circ c^ {-1}$
as in Figure \ref{fig:uniquepants}.

Then $\G(\delta)$ is either a simple closed geodesic or another
cone point, and in either case $(p,\gamma,\G(\gamma))$ is a pair
of pants. Using exactly the same technique with two cone points
and a simple geodesic path between them, or with two simple
closed geodesics and a perpendicular simple geodesic path between
them, we get pairs of pants that are uniquely determined. This
construction is essential for the following theorem.

\begin{thm}\label{thm:collartheorem}

Let $M$ be an admissible cone-surface of genus $g$ with $n$ cone
points $p_1,\ldots,p_n$ with cone angles
$2\varphi_1,\ldots,2\varphi_n$. Let $2\varphi$ be the largest
cone angle. Let $\gamma_1,\ldots,\gamma_m$ be disjoint simple
closed geodesics on $M$. Then the following hold.
\begin{enumerate}
\item $m\leq 3g-3 +n$.\\

\item There exist simple closed geodesics
$\gamma_{m+1},\ldots,\gamma_{3g-3+n}$ which together with
$\gamma_1,\ldots,\gamma_m$ form a partition of $M$.\\

\item The collars
$$
\CC(\gamma_k)=\{x\in M \mid d(x,\gamma_k) \leq w_k=\arcsinh (\cos
\varphi /\sinh \frac{\gamma_k}{2})\}
$$

\noindent and

$$
\CC(p_l)=\{x\in M \mid d(x,p_l) \leq v_l=\arccosh(1 / \sin
\varphi_l)\}
$$

\noindent are pairwise disjoint for $k = 1,\ldots,3g-3+n$ and
$l=1,\ldots,n$.\\

\item

Each $\CC(\gamma_k)$ is isometric to the cylinder
$[-w_k,w_k]\times \Sp^1$ with the Riemannian metric
$ds^2=d\rho^2+\ell^2(\gamma_k) \cosh^2\rho dt^2$.\\

\noindent Each $\CC(p_l)$ is isometric to a hyperbolic cone
$[0,v_l]\times \Sp^1$ with the Riemannian metric
$ds^2=d\rho^2+\frac{\varphi_l^2}{\pi^2}\sinh^2\rho dt^2$.
\end{enumerate}

\end{thm}

\begin{proof}

The first two points are equivalent to the problem of counting the
number of geodesics in a partition for a surface of signature
$(g,n)$, and showing that any collection of pairwise disjoint simple
closed geodesics can be completed to form a partition.  These
questions are not new and the proofs are known (e.g. \cite{bubook}).

The first step in proving the theorem is to show that a cone
point $p_l$ is at a distance of at least $\arccosh(1 / \sin
\varphi_l)$ from all simple closed geodesics. Let $\gamma$ be a
simple closed geodesic on $M$. Take a geodesic path $c$ that
realizes the distance between $p_l$ and $\gamma$ (i.e. such that
$\ell(c)=d(p_l,\gamma)$). Then take the unique pair of pants
obtained from $p_l,\gamma$ and $\G(\gamma\circ c \circ \eta \circ
c^{-1})$ as discussed previously. Either this pair of pants is a
$V$-piece or a joker's hat. In both cases extract the
trirectangle (a hyperbolic quadrilateral with three right angles)
as in Figures \ref{fig:collar1} and \ref{fig:collar2}.

\begin{figure}
\center
\epsfxsize=50mm 
\epsfbox{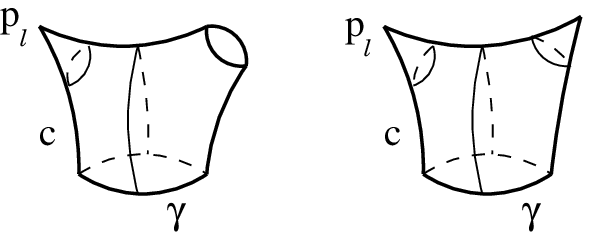}
\caption{}
\label{fig:collar1}
\end{figure}

\begin{figure}
\center
\epsfxsize=20mm 
\epsfbox{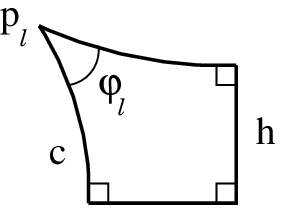}
\caption{}
\label{fig:collar2}
\end{figure}

From the formulas for the trirectangle we have
$$
\cosh c= \frac{\cosh h}{\sin \varphi_l}
$$
where $h$ is as in Figure \ref{fig:collar2}. It follows that
$c>\arccosh(1 / \sin
\varphi_l)$.

Let $p_l$ and $p_{l'}$ be two distinct cone points. Let $b$ be a
path that realizes distance between them. Let $\gamma$ be a
simple closed geodesic that crosses $b$. It is clear that
$\ell(b)$ is necessarily greater than or equal to the distance
from $p_l$ to $\gamma$ added to the distance from $p_{l'}$ to
$\gamma$. From what precedes we have:
$$
\ell(b)\geq
d(p_l,\gamma)+d(p_{l'},\gamma)>\arccosh\left(\frac{1}{\sin
\varphi_l}\right) + \arccosh\left(\frac{1}{\sin \varphi_ {l'}}
\right).
$$
It follows that the distance sets $\CC(p_l)$ and $\CC(p_ {l'})$
are disjoint.

Let $\gamma_k$ and $\gamma_{k'}$ be two disjoint simple closed
geodesics. Let $c$ be a path that realizes distance between them.
The unique pair of pants determined by $c$, $\gamma_k$ and
$\gamma_{k'}$ is either a $Y$-piece or a $V$-piece. In the first
case it follows that $\CC(\gamma_k)$ and $\CC(\gamma_{k'})$ are
disjoint sets from the collar theorem on Riemann surfaces. In the
latter case, consider Figures \ref{fig:collar3} and \ref{fig:collar4}.  

\begin{figure}
\center
\epsfxsize=25mm 
\epsfbox{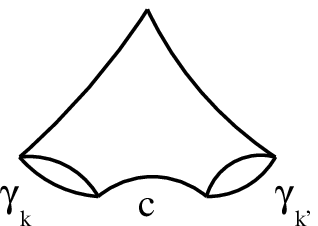}
\caption{}
\label{fig:collar3}
\end{figure}

\begin{figure}
\center
\epsfxsize=35mm 
\epsfbox{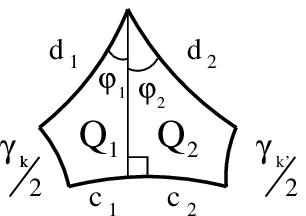}
\caption{}
\label{fig:collar4}
\end{figure}

Notice that only half of the collar around a given boundary geodesic
is contained in a pair of pants.
Also note that both angles $\varphi_1$ and $\varphi_2$ are strictly
inferior to $\varphi$. From the trirectangle formula applied to
$Q_1$ we obtain
$$
\sinh c_1 = \frac{\cos \varphi_1}{\sinh\frac{\gamma_k}{2}},
$$

\noindent and analogously for $Q_2$. From this we obtain that

$$
c_1 = \arcsinh\left(\frac{\cos
\varphi_1}{\sinh\frac{\gamma_k}{2}}\right),
$$

\noindent and
$$
c_2 = \arcsinh\left(\frac{\cos
\varphi_2}{\sinh\frac{\gamma_{k'}}{2}}\right). $$
It follows that the distance sets $\CC(\gamma_k)$ and
$\CC(\gamma_{k'})$ are disjoint, because each collar is the union
of two such half-collars.

It remains to prove that for arbitrary $p_l$ and $\gamma_k$ the
collars are disjoint. Let $c$ be a geodesic path that realizes
the distance between them. The collars around both the cone point
and the geodesic have widths which depend only on their angle or
length. It thus suffices to consider the case where the cone point
and the geodesic would be as close as possible. This case would
occur if the pair of pants resulting from $p_l$, $\gamma_k$ and
$c$ were a joker's hat with both cone angles equal to $2\varphi$.
(This is a direct consequence of hyperbolic geometry.) In this
case Figure \ref{fig:collar5} would apply.

\begin{figure}
\center
\epsfxsize=20mm 
\epsfbox{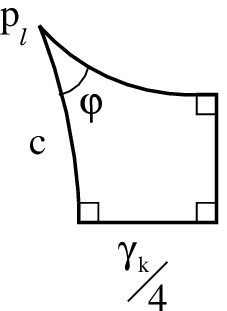}
\caption{}
\label{fig:collar5}
\end{figure}

From a trirectangle formula one obtains
$$
\sinh c = \coth \frac{\gamma_k}{4} \cot \varphi.
$$
This in turn can be expressed as
$$\sinh c = \frac{\cosh
\frac{\gamma_k}{2}+1}{\sinh\frac{\gamma_k}{2}} \cot \varphi.
$$

Let us compare this with the value obtained by calculating
$\sinh(w_k+v_l)$. By calculation one obtains

$$
\sinh(w_k+v_l)=\frac{1+\sqrt{\cos^2\varphi+\sinh^2\frac
{\gamma_k}{2}}}{\sinh\frac{\gamma_k}{2}} \cdot \cot \varphi.
$$
By comparison and because $\varphi < \frac{\pi}{2}$ we obtain that
$c > w_k + v_l$. This implies that the distance sets
$\CC(\gamma_k)$ and $\CC(p_l)$ are disjoint.

The last point of the theorem is obtained as in the classical
theory of Riemann surfaces for $\CC(\gamma_k)$ (see \cite{bubook}),
and was shown by Dianu \cite{diphd} for $\CC(p_l)$.
\end{proof}

\begin{rem}
The values for the collars are in general optimal in the following sense. The
collar around a simple closed geodesic $\gamma$ can be seen as a
distance set with the following property: if another simple closed
geodesic enters $\CC(\gamma)$, then it necessarily intersects
$\gamma$. Suppose $\gamma$ is a non-separating simple closed geodesic on an
admissible cone-surface with $(g,n)\neq (1,1)$. (If $(g,n)=(1,1)$ then all simple
closed geodesics intersect.) Now
replacing $w_\gamma=\arcsinh (\cos \varphi /\sinh \frac{\gamma}{2})$ by
$w=w_\gamma+\varepsilon$ with $\varepsilon > 0$ would be fatal to
this intersection property. In fact, for any $\varepsilon > 0$ there are an infinity of simple 
closed geodesics that intersect the enlarged collar but {\it not} $\gamma$. To prove
this, take a $V$-piece containing $\gamma$ and $p$ where the cone
angle at $p$ is exactly $2\varphi$. The other boundary geodesic
$\gamma'$ of the $V$-piece can be chosen as long as desired (see
\cite{pa041}). From the formulas obtained in the proof, it follows
that $\gamma'$ can be arbitrarily close to the collar of $\gamma$.
In the case where $\gamma$ is a separating geodesic, the optimal collar
width will depend on how the cone angles are distributed on both sides of
$\gamma$; using the same techniques, one could find an optimal collar in 
each of the various cases. 

In an analogous fashion, one can show that the bound for
collars around cone points is also sharp, provided one can include the cone
point in a $V$-piece with two distinct simple closed geodesics. This excludes
only the cases where $(g,n)=(1,1)$ and $(0,4)$. In the latter case, 
the optimal collar width depends on all four cone angles. 
Now if the surface is a torus $T$ with a single cone point $p$ with cone angle $2\varphi$,
by cutting along an interior simple closed geodesic of $T$, say
$\gamma$, one obtains a $V$-piece with two boundary geodesics of
equal length. The process of lengthening only one of the two
geodesics is no longer possible and this extra rigidity implies
that the sharp constant is in fact $\arccosh (1/\sin
\frac{\varphi}{2})$ (cf. \cite{pa041}).
\end{rem}

An example of the utility of the collar theorem is the following
natural corollary.

\begin{cor}
Let $\gamma$ and $\delta$ be closed geodesics on $M$ which
intersect each other transversally, and assume that $\gamma$ is
simple. Then
$$
\sinh \frac{\ell(\gamma)}{2} \sinh \frac{\ell(\delta)}{2} > \cos
\varphi.
$$
\end{cor}

\section{Partitions} \label{Sect:S4}

We now proceed to prove Bers' theorem for compact admissible
cone-surfaces.

\begin{thm}\label{thm:bersthm}
Let $M$ be an admissible cone-surface of genus $g$ with
$n$ cone points. Then there exists a partition $\PP$ of $M$ such
that every geodesic in $\PP$ has length less than a constant
$L_{g,n}$.
\end{thm}

\begin{proof}
Let $p_1, \ldots, p_n$ be the cone points on $M$ with
corresponding cone angles $2\varphi_1, \ldots, 2\varphi_n$. We
define
\begin{displaymath}
Z_i (r_i) = \{ x \in M | \mbox{\ dist} (x, p_i) \leq r_i \},
\end{displaymath}
for $i = 1, \ldots, n$. We denote the boundary of $Z_i (r_i) $ by
$\beta_i$.  Each neighborhood $Z_i(r_i)$ admits polar coordinates, and we have
\begin{displaymath}
\mbox{Area } Z_i (r_i) = 2 \varphi_i (\cosh r_i - 1).
\end{displaymath}
For $i = 1, \ldots, n$, we also know that
\begin{displaymath}
\mbox{Area } Z_i(r_i) < \mbox{Area } M.
\end{displaymath}
This implies
\begin{displaymath}
r_i < \arccosh \left( 1 + \frac{\mbox{Area } M}{2 \varphi_i} \right).
\end{displaymath}
Now
\begin{displaymath}
\ell(\beta_i) = 2 \varphi_i \sinh r_i,
\end{displaymath}
so
\begin{eqnarray*}
\ell (\beta_i) & < & 2 \varphi_i \sinh \left(\arccosh
\left(1 + \frac{\mbox{Area } M}{2 \varphi_i}\right)\right) \\
& = & 2 \varphi_i \sqrt{ \left(1 + \frac{\mbox{Area } M}{2
      \varphi_i}\right)^2- 1} \\
& < & \mbox{Area } M + 2 \varphi_i \\
& < & 2 \pi (2g-2+n).
\end{eqnarray*}
Thus the length of the boundary of the neighborhood $Z_i (r_i)$
about $p_i$, $i = 1, \ldots, n$, is bounded above by $2 \pi
(2g-2+n)$.

Starting from $r_i = 0$, we let $r_i$ grow
continuously (for all $i$ simultaneously) until one of the
following cases occurs: 

\begin{enumerate}
\item $\beta_j$ ceases to be simple for some $j \in \{1,
\ldots, n\}$;
\item $\beta_j$ meets $\beta_k$ for some $j \neq k$.
\end{enumerate}
Once one of these cases occurs, we fix all the $r_i$.  We will
consider each case in turn. \\

\noindent \textbf{Case 1}.  We view $\beta_j$ as the composition of
two curves
$\tau_1$ and $\tau_2$, which both have initial and final points
at the self-intersection of $\beta_j$. If $\beta_j$ has
multiple self-intersections, then we choose one and let $\tau_1$
and $\tau_2$ be as in Figure \ref{fig:necklace}.

\begin{figure}
\center
\epsfxsize=50mm 
\epsfbox{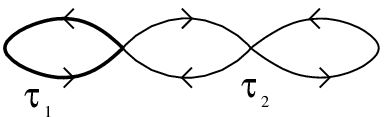}
\caption{}
\label{fig:necklace}
\end{figure}

We allow $\tau_1$ and $\tau_2$ to slide in their free homotopy
classes on $M \setminus \Sigma$ until they reach their minima.
We have
\begin{displaymath}
\ell (\G (\tau_1)) < \ell ( \tau_1) < \ell(\beta_j) < 2 \pi
(2g-2+n),
\end{displaymath}
and similarly for $\ell (\G (\tau_2))$.  Note that it is possible that
one of $\G (\tau_1)$ or $\G (\tau_2)$ is a cone point (but not both,
as $M$ has signature $(g,n) \geq (0,4)$), in which case
the nontrivial geodesic bounds a joker's hat.  
We cut $M$ open along those $\G (\tau_i), \ i=1,2,$ which are not cone
points, and remove any
resulting $V$-pieces or joker's hats.  Let $M^1$ be the
(possibly empty) remaining connected component.  
Then 
$$
\ell(\partial M^1) < 4 \pi (2g-2+n).
$$ \\

\noindent \textbf{Case 2}.  Suppose that two boundary
curves $\beta_j$ and $\beta_k$ meet, and that both
$\beta_j$ and $\beta_k$ are simple (if more than two simple boundary curves
meet, choose two). Consider the
curve $\tau$ obtained by first traversing $\beta_j$ and then
traversing $\beta_k$, where the initial and final points of
both of these curves are at their intersection point. Note that
$\tau$ is homotopic to a simple closed curve and is not homotopic to a
cone point, as $M$ has signature $(g,n) \neq (0,3)$. 
Thus, $\tau$ is homotopic to a unique simple closed
geodesic. We have
\begin{displaymath}
\ell (\G (\tau)) < \ell (\tau) = \ell (\beta_j) + \ell
(\beta_k) < 4 \pi (2g-2+n).
\end{displaymath}

Cutting $M$ open along this geodesic yields at least one
joker's hat.  Let $M^1$ denote the (cone-)surface obtained by cutting
$M$ open in this way and removing any joker's hats that result.  Then
$$
\ell(\partial M^1) < 4 \pi (2g-2+n).
$$  

We now restart the process; that is, we send out collars $Z_i(r_i)$
from the remaining cone points $p_i$ on $M^1$ and let $r_i$ grow
continuously from $r_i=0$ until one of the following situations occurs:

\begin{enumerate}
\item $\beta_j$ ceases to be simple for some $j \in \{1,\ldots, n\}$;
\item $\beta_j$ meets $\beta_k$ for some $j \neq k$;
\item $\beta_j$ meets a boundary geodesic $\gamma_i$ on $M^1$. 
\end{enumerate}
Once one of these cases occurs, we fix all the $r_i$.  Cases 1 and 2
are as above, and if we let $M^2$ be the (cone-)surface which results
from cutting $M^1$ open along the new geodesics we find and removing
all Y-pieces, V-pieces and joker's hats, then 
$$
\ell( \partial M^2) < 8 \pi (2g-2+n).
$$  \\

\noindent \textbf{Case 3}.  Consider the curve $\tau$ obtained by first
traversing $\beta_j$ and then traversing $\gamma_i$, 
where the initial and final points of
both of these curves are at their intersection point. Note that
$\tau$ is homotopic to a simple curve and is not homotopic to a
cone point; if it were, then $p_j$ would live on a joker's hat that
would have been removed at the previous step.  
Thus, $\tau$ is homotopic to a unique simple closed
geodesic. We have
\begin{displaymath}
\ell (\G (\tau)) < \ell (\tau) = \ell (\beta_j) + \ell
(\gamma_i) < 2 \pi (2g-2+n) + 4 \pi (2g-2+n) = 6 \pi (2g-2+n).
\end{displaymath}
Cutting $M^1$ open along $\G (\tau)$ yields at least one V-piece.  Let
$M^2$ be as defined above; then
$$
\ell (\partial M^2) < \ell(\partial M^1) + 2 \pi (2g-2+n) < 6 \pi (2g-2+n),
$$
as $\ell(\G(\tau)) < \ell(\gamma_i) + 2 \pi (2g-2+n)$ and $\gamma_i
\notin \partial M^2$.

We repeat the above process until all of the cone points on $M$
have been removed on V-pieces or joker's hats.  Note that at each
step, after cutting our cone-surface open along the geodesics we
find and removing any Y-pieces, V-pieces and joker's hats, the length of the
boundary of the resulting (cone-)surface increases by at most $4 \pi
(2g-2+n)$.  To remove all $n$ cone points requires $m \leq 2n$
geodesics and $\mu \leq m$ steps; thus, we have found 
$\gamma_1, \ldots, \gamma_m$ such that
$$
\ell (\partial M^j) < 4 \pi j (2g-2+n), \mbox{\ \ \ \ } j=1, \ldots, \mu
$$
and
$$
\ell (\gamma_k) < 4 \pi k (2g-2+n), \mbox{\ \ \ \ } k=1, \ldots,m. 
$$

To obtain the remaining geodesics in our decomposition, we
proceed by induction. That is, we find a suitable simple closed
geodesic in the interior of $M^{\mu}$ which is not homotopic to a
boundary component, cut $M^{\mu}$ open along this geodesic, remove any
$Y$-pieces, and let $M^{\mu+1}$ be the resulting surface. To find such
a geodesic, we create tubular neighborhoods around all boundary
geodesics and let the widths grow until a critical case occurs;
the area arguments are analogous to those in the induction for
compact Riemann surfaces of genus $g \geq 2$ (see \cite {bubook}). 
\end{proof}

\begin{rem} 
The proof gives an explicit bound for the 
length of each partitioning geodesic:
\begin{displaymath}
\ell(\gamma_k) < 4 \pi k (2g-2+n),
\end{displaymath}
where $\gamma_k$ is the $k$th geodesic in a partition of $M$. For
$L_{g,n}$ we have thus proved the following bound:
$$
L_{g,n}< 4 \pi (3g-3+n) (2g-2+n).
$$
\end{rem}

\bibliographystyle{plain}
\bibliography{collarandbers.bib}

\end{document}